# Penalized estimate of the number of states in Gaussian linear AR with Markov regime

## Ricardo Ríos

*Departamento de Matemáticas*
*Escuela de Matemáticas*
*Universidad Central de Venezuela*
*Caracas-Venezuela*
*e-mail:* ricardo.rios@ciens.ucv.ve

## Luis-Angel Rodríguez

*Departamento de Matemáticas*
*Facultad de Ciencias y Tecnología*
*Universidad de Carabobo*
*Valencia-Venezuela*
*e-mail:* larodri@uc.edu.ve

**Abstract:** We deal with the estimation of the regime number in a linear Gaussian autoregressive process with a Markov regime (AR-MR). The problem of estimating the number of regimes in this type of series is that of determining the number of states in the hidden Markov chain controlling the process. We propose a method based on penalized maximum likelihood estimation and establish its strong consistency (almost sure) without assuming previous bounds on the number of states.

## 1. Introduction

Our aim in this article is to establish consistency criteria for the method of penalized maximum likelihood estimation for the number of states in a hidden Markov chain in an AR-MR process. We show strong consistency for an estimator of the number of states in autoregressive process with Markov regime when the regression functions are linear and the noise is Gaussian.

Autoregressive processes with Markov regime can be looked at as a combination of hidden Markov models (HMM) with threshold regression models. These have been introduced in an econometric context by Goldfeld and Quandt [12] and they have become quite popular in the literature ever since Hamilton [13] employed them in the analysis of the gross internal product of the USA for two regimes: one of contraction and another of expansion.

When the number of states in a hidden Markov chain is known *a priori* the estimation problems can be solved, in principle, through the use of techniques based on maximum likelihood estimation (see, McDonald and Zucchini [17] and Cappe *et al.* [2]). But in many applications, a key problem is to determine the number of states in a way such that the data is adequately described while at the same time a compromise is maintained between fitness and possibility of generalizing the model. The problem of estimating the hidden Markov chain in AR-MR is a typical example of a nested-family of models: models with $m$ parameters can be seen also as models with $m+1$ parameters. Thus the problem of model selection is essentially that of determining the smallest model that contains the distribution capable of generating the data. In many instances, the estimation of the model will depend on how identifiability affects the model and not on the specification of the correct model.

A first approximation to determine the dimension of the model is a statistical test based on the likelihood ratio (see Dachuna and Duflo [5], p. 227). For the estimation of the number of state in hidden Markov chain, the likelihood ratio test fails because regularity assumptions do not hold. In particular, the model is not identifiable, as some parameters do not show up under the null hypothesis and the information matrix is singular. As a result the asymptotic distribution of the likelihood ratio is not $\chi^2$. As an alternative, one can construct generalized tests for the likelihood ratio that would hold under non-standard conditions. For the problem of the determination of the number of states in AR-MR, Hansen [14] has proposed a test that works with loss of identifiability but in order to implement it one needs to calculate *p*-values in an approximate way; this leads to computationally heavy calculations which produce approximate *p*-values which underestimate the real ones. Garcia [8] has advanced more attractive computational alternatives which lack, however, the technical rigor present in Hansen's approach.

For HMM models the likelihood ratio test is not bounded. Gassiat and Keribin have studied it [11] and have shown that it diverges to infinity. The rate of growth of the likelihood ratio as the parameters increase is related to the complexity of the model. This brings us to consider penalized estimators of the likelihood function that compensate the lack of likeness between models with different dimensions. The specification of small penalties depends on how the divergence rate at infinitum of the likelihood ratio is determined. But as far as we know this is still an open problem for HMM models where the data belongs to infinite sets.

In general, criteria for penalized likelihood are obtained through approximations to Kullback-Leibler divergence. Among others, we find the very popular information criteria of Akaike (AIC) and the Bayesian one (BIC). These have been used by several authors in applications of the HMM models, however, as is mentioned by McDonald and Zucchini [17], these authors have made no reference as to their validity.

We shall distinguish two cases, regarding whether or not the observed variables are in an infinite set. For the case of the HMM model with data belonging to a finite set much work has been done starting with Finesso's presentation of the problem [7] where he establishes the strong consistency for the penalized estimator of the number of states assuming that the actual number of states belongs to bounded set of integers. Liu y Narayan [16], assuming this restriction introduce a strongly-consistent estimator based on statistical mixtures of the Krischevsky-Trofimov type as this allows to normalize the likelihood so as to control the growth of likelihood when the number of states is increased. In studies dealing with the efficiency of their estimator they show specifically that the probability for underestimating decreases at an exponential rate with the sample size, whereas the probability for overestimating does not exceed a third degree polynomial of the size of the sample. Based on this former work, Gassiat and Boucheron [10] have introduced considerable advances: they proved the strong consistency of the penalized estimator without assuming *a priori* upper bounds for the number of states; in addition, they showed that the probabilities for underestimating as well as for overestimating fall at an exponential rate with sample size. For AR-MR processes with observations belonging to a finite set the techniques introduced by Gassiat and Boucheron were further used by Chambaz and Matias [4] to simultaneously show the consistency of the number of states of the hidden chain and the memory of the observed process.

For the non-finite case in HMM models Rydén [19] have shown consistency for a penalized likelihood estimator which in the limit does nor underestimate the number of states. Dortet-Bernadet [6] have shown that under certain regularity conditions the Rydén estimator is indeed consistent. Gassiat [9] studying a penalized estimator of marginal likelihood concludes that there is consistency in probability with the actual number of states. This technique is extended by Olteanu and Rynkiewicz [18] in order to select the number of regression functions in processes where the regime is controlled by an independent sequence. In this very same work, the authors indicate that the penalized marginal likelihood criterion cannot be directly applied to AR-MR. Smith *et al.* [21] have advanced a new information criterion in order to be able to approximate the Kullback-Leibler divergence and to select the numbers of states and the variables in AR-MR. This criterion imposes a penalty that reduces state number overestimation. Following the work on finite alphabets in Ref. [7, 16, 10], Chambaz *et al.* [3] have shown strong consistency for penalized and Bayesian estimators of the number of states in HMM and observations belonging to infinite (discrete and continuum) sets; they have worked with conditionally Poisson and Gaussian distribution. As in the previous works, no *a priori* bounds are assumed for the number of states.

Following Chambaz *et al.* [3] we prove a mixture-type inequality (see Section 2.1) that allows us to normalize the likelihood and in addition we also prove in Section 3, without assuming *a priori* bounds on the actual state number of the hidden Markov chain, that the penalized estimator underestimates. In order to show that the penalized estimator does not overestimate the number of states, we use an approach that works well for nesting models and which is based on the equicontinuity of the likelihood function. We would like to point out that our results are obtained for the linear case and that they can be easily generalized to the nonlinear case if we assume that a sublinearity hypothesis such as the one required by Yao y Attali [22] holds, albeit retaining the assumption of Gaussian-like behavior.

## 2. Definitions and introductory comments

A linear autoregressive process with Markov regime (AR-MR) is defined by:

$$Y_n = \alpha_{X_n} Y_{n-1} + b_{X_n} + \sigma_{X_n} e_n \tag{2.1}$$

where $\{e_n\}$ are i.i.d. random variables, $\sigma_i^2$ is the variance of the model in each regime and $\sigma^2 = (\sigma_1^2, \ldots, \sigma_m^2)$. The sequence $\{X_n\}$ is a homogeneous Markov chain with state space $\{1, \ldots, m\}$. We denote by $A$ its transition matrix $A = [a_{ij}]$. For each $1 \leq i \leq m$ we have $\theta_i = (b_i, \alpha_i)^t$ and

$$\theta = \begin{pmatrix} b_1 & b_2 & \cdots & b_m \\ \alpha_1 & \alpha_2 & \cdots & \alpha_m \end{pmatrix}.$$

We assume that:

**S1** The Markov chain $\{X_n\}$ is recurrent positive. Hence, it can have an invariant distribution that we denote by $\lambda = (\lambda_1, \ldots, \lambda_m)$.
**S2** $Y_0$, the Markov chain $\{X_n\}$ and the sequence $\{e_n\}$ are mutually independent.
**S3** The $e_n$ has Gaussian distribution $\mathcal{N}(0, 1)$.
**S4** $\mathbb{E}_\lambda(\log \alpha) = \sum_{i=1}^m \lambda_i \log(\alpha_i) < 0$ (stability condition).
**S5** The parameter $\theta_i$ belongs to the compact subset $\Theta_i \subset \mathbb{R}^2$.
**S6** For each $1 \leq i \leq m$, $\sigma_i^2 \in [c, d]$, $c > 0$.

The parameter space is the set

$$\Psi_m = \left\{ \psi = (\theta, \sigma^2, A) : \theta \in \bigotimes_{i=1}^m \Theta_i,\ \sigma^2 \in [c,d]^m,\ \sum_{j=1}^m a_{ij} = 1 \right\}.$$

Notations

- $V_1^n$ stands for random vector $(V_1, \ldots, V_n)^t$ and $v_1^n = (v_1, \ldots, v_n)^t$ for any realization.
- The symbol $\mathbb{1}_B(x)$ denotes the function which assigns the value 1 if $x \in B$ and 0 elsewhere.
- Distributions and densities are denote by $p$.

For each $1 \leq i \leq m$,

- Let $n_i = \sum_{k=1}^n \mathbb{1}_j(x_k)$ be the number of visits of a realization of the Markov chain $\{X_n\}$ to state $i$ in the first $n$ steps. $n_{ij} = \sum_{k=1}^{n-1} \mathbb{1}_{i,j}(x_{k-1}, x_k)$ is the number of transitions from $i$ to $j$ in $n$ steps.
- Let $I_i := \{k \leq n : X_k = i\} = \{k_{1_i}, \ldots, k_{n_i}\}$.

- Let

$$\begin{aligned}
\mathbf{Y}_{I_i} &:= (Y_{k_{1_i}}, \ldots, Y_{k_{n_i}})^t \\
\mathbf{Y}_{I_{i-1}} &:= (Y_{k_{1_i}-1}, \ldots, Y_{k_{n_i}-1})^t \\
\mathbf{E}_i &:= \{e_{k_{1_i}}, \ldots, e_{k_{n_i}}\}.
\end{aligned}$$

- The symbol $\mathbf{1}_i$ denotes a $n_i$-dimensional column vector with 1 in all of its positions and $\mathbf{W}_i = [\mathbf{1}_i, \mathbf{Y}_{I_i}]$.

The process $\{Y_n\}$ in general is not a Markov chain but the associated process $\{(Y_n, X_n)\}$ is a Markov chain with state space $\mathbb{R} \times \{1, \ldots, m\}$. In what follows we introduce some properties– to be used throughout this work– related to the likelihood function for the present model.

## 3. The likelihood function

We consider the conditional distribution $p_\psi(Y_1^n|Y_0 = y_0)$ as the likelihood function for a set of observations $Y_0^n = y_0^n$ and parameter $\psi$. Because of the total probability rule the total, likelihood function for the model is given by:

$$\begin{aligned}
&p_\psi(Y_1^n|Y_0 = y_0) \\
&= \sum_{x_1^n} p_\psi(Y_1^n, x_1^n | Y_0 = y_0) = \sum_{x_1^n} p_{\theta,\sigma^2}(Y_1^n|Y_0 = y_0, x_1^n) p_A(x_1^n). \quad (3.1)
\end{aligned}$$

Using our above notation we may represent the AR-MR process defined by Eq. (2.1) by means of its $m$ linear models, for each $1 \leq i \leq m$

$$\mathbf{Y}_{I_i} = \mathbf{W}_i \theta_i + \sigma_i \mathbf{E}_i \quad \forall i \leq m. \quad (3.2)$$

Thus, the distribution of $Y_0^n$ conditional to $x_1^n$ is written as

$$p_\psi(Y_1^n|Y_0 = y_0, x_1^n) = \prod_{i=1}^m \frac{1}{(\sqrt{2\pi\sigma_i^2})^{n_i}} \exp\left(-\frac{1}{2\sigma_i^2}(\mathbf{Y}_{I_i} - \mathbf{W}_i\theta_i)^t(\mathbf{Y}_{I_i} - \mathbf{W}_i\theta_i)\right).$$

We assume that prior distribution $p(\psi)$ on $\Psi$ satisfies

$$p(\psi) = p(A)p(\theta|\sigma^2)p(\sigma^2) = \prod_{i=1}^m p(A_i)p(\theta_i|\sigma_i^2)p(\sigma_i^2),$$

where $A_i$ denotes the $i$-th row of $A$. Due to (3.2) we will consider the prior distribution for $(\theta, \sigma^2)$ belonging to a Gaussian-Gamma family (see Broemiling [1], §1, page. 3), means for each $i = 1 \ldots, m$,

**H1** $\theta_1, \ldots, \theta_m$ are independent with

$$\theta_i \sim \mathcal{N}(\theta_i|0, \sigma_i^2\tau^2\mathbb{I}) = \frac{1}{2\pi\sigma_i^2\tau^2} \exp\left(-\frac{1}{2\sigma_i^2\tau^2}\theta_i^t\theta_i\right)$$

**H2** $\sigma_1^2, \ldots, \sigma_m^2$ are independent with Inverse-Gamma distribution

$$\sigma_i^2 \sim \mathcal{IG}(v_0/2, u_0/2) = \frac{\left(\frac{u_0}{2}\right)^{v_0/2}}{\Gamma(v_0/2)}(\sigma_i^2)^{-(\frac{v_0}{2}+1)}e^{-\frac{u_0}{2\sigma_i^2}}.$$

**H3** $A_1, \ldots, A_m$ are independent with $A_i \sim \mathcal{D}(e_i)$, where $\mathcal{D}$ denotes a Dirichlet density the parameters vector $(1/2, \ldots, 1/2)$,

$$\mathcal{D}(e_i) = \frac{\Gamma(m/2)}{\Gamma(1/2)^m}\prod_{j=1}^{m}a_{ij}^{-1/2}.$$

The related mixture statistic is defined by

$$q_m(Y_1^n) = \int_\Psi p_\psi(Y_1^n|Y_0 = y_0)p(\psi)d\psi.$$

The main results of this section is the comparison between the likelihood function and the mixture statistics.

Under the assumptions (**S1-S6**) and (**H1-H3**) described before we have the following theorem.

**Theorem 3.1.** *For each $m \geq 1$ and the prior distribution $p(\psi)$ satisfies the inequality*

$$\log \frac{p_\psi(Y_1^n|Y_0 = y_0)}{q_m(Y_1^n)}$$
$$\leq \frac{m(m+1)}{2}\log(n) + c_m(n) + d(n) + \frac{nm}{2}\log\frac{\mathbf{Y}_{I_k}^T\mathbf{P}_k\mathbf{Y}_{I_k}}{\mathbf{Y}_{I_k}^t\mathbf{B}_k\mathbf{Y}_{I_k}} + e_m(n),$$

*where*

$$\frac{\mathbf{Y}_{I_k}^T\mathbf{P}_k\mathbf{Y}_{I_k}}{\mathbf{Y}_{I_k}^t\mathbf{B}_k\mathbf{Y}_{I_k}} = \max_{i=1,\ldots,m}\frac{\mathbf{Y}_{I_i}^T\mathbf{P}_i\mathbf{Y}_{I_i}}{\mathbf{Y}_i^t\mathbf{B}_i\mathbf{Y}_{I_i}}$$

*and for each $n \geq 4$,*

$$c_m(n) = \max\left\{0, \log m - m\left(\log\frac{\Gamma(m/2)}{\Gamma(1/2)} - \frac{m(m-1)}{4n} + \frac{1}{12n}\right)\right\},$$
$$d(n) = \frac{n}{2} + \frac{1}{2}\log\left(\frac{n}{2}\right),$$
$$e_m(n) = \max\left\{0, \frac{m}{2}\log\left(\frac{1}{n^2} + \frac{\tau^4}{m}\sum_{i=1}^{m}(\lambda_i\sigma_i)^2\right) - \frac{m\log(2\pi)}{2}\right\}$$
$$\mathbf{P}_i = \mathbb{I} - \mathbf{W}_i\mathbf{M}_i\mathbf{W}_i^T$$
$$\mathbf{M}_i = (\mathbf{W}_i^T\mathbf{W}_i + \tau^{-2}\mathbb{I})^{-1}$$
$$\mathbf{B_i} = \mathbb{I} - \mathbf{W}_i(\mathbf{W}_i^T\mathbf{W}_i)^{-1}\mathbf{W}_i^T.$$

## 4. Penalized estimation of the number of states

The purpose of this Section is to advance an estimation method based on penalized maximum likelihood in order to select the number of states $m$ of a hidden Markov chain $\{X_n\}$. For every integer $m \geq 1$, we consider the sets $\Psi_m$ and $\mathcal{M} = \bigcup_{m\geq 1} \Psi_m$ the family for all models, (with convention $\Psi_0 = \emptyset$). We define the number of states $m_0$ through the property

$$p_{\psi_{m_0}} \in \{p_\psi : \psi \in \Psi_{m_0}\} \cap \{p_\psi : \psi \in \Psi_{m_0-1}\}^c. \tag{4.1}$$

Remark: (Identifiability) We assume that for the true model $\Psi_{m_0}$ the vector components $\{(\alpha_i, b_i, \sigma_i)\}_{i=1}^{m_0}$ are different; thus, for every $n$, there exists a point $Y_{n-1} \in \mathbb{R}$ such that $\{(\alpha_i Y_{n-1} + b_i, \sigma_i)\}_{i=1}^{m_0}$ are different. Therefore, in agreement with Remark 2.10 of Krishnamurthy and Yin [15] the model is identifiable in the following sense: If $K$ stands for the Kullback-Leibler divergence $K(\psi, \psi_{m_0}) = 0$ then, $\psi = \psi_{m_0}$. As a result, identifiability implies that $m_0$ – defined by Eq. (4.1) – is unique.

Let $pen(n,m)$ be a penalty term which is given by a positive function with increasing values of $n$ and $m$. We define the estimator for penalized maximum likelihood as (PML) for $m_0$ as,

$$\widehat{m}(n) = \underset{m\geq 1}{\operatorname{argmin}} \Big\{ - \sup_{\psi \in \Psi} \log p_\psi(Y_1^n | Y_0 = y_0) + pen(n,m) \Big\}. \tag{4.2}$$

We say that $\widehat{m}(n)$ overestimates the number of states $m_0$ if $\widehat{m}(n) > m_0$ and that it underestimates the number of states if $\widehat{m}(n) < m_0$.

In the following theorem we prove that the estimator PML for $m_0$, overestimates the number of states.

**Theorem 4.1.** *Assume (S1-S6) and that* $\lim_{n\to\infty} \frac{pen(n,m)}{n} = 0 \ \forall \ m$. *Then*

$$\widehat{m}(n) \geq m_0. \ a.s.$$

In order to prove this Theorem, the following two Lemmas are necessary:

**Lemma 4.1** (Finesso [7]). *Assume (S1-S6) the set of functions $f_n(\psi) = \frac{1}{n} \log p_\psi(Y_1^n | Y_0 = y_0)$ is an equicontinuos sequence a.s.-$\mathbb{P}_{\psi_0}$.*

The following result is a usual one in the context of order selection for a nested family of models, see [2], §15, p 577–578. For HMM models, similar results are given, for example, in [10, 3].

**Lemma 4.2.** *Assume (S1-S6) we have:*

1. *For each $m \geq 1$, $\psi, \psi_0 \in \Psi_m$ there exist $K(\psi, \psi_0) < \infty$ such that:*

$$\lim_{n\to\infty} [\log p_{\psi_0}(Y_1^n | Y_0 = y_0) - \log p_\psi(Y_1^n | Y_0 = y_0)] = K(\psi, \psi_0).$$

2. *For each $\psi \in \Psi_{m_0} \cap \Psi_{m_0-1}^c$,*

$$\min_{m < m_0} \inf_{\psi \in \Psi_m} K(\psi_{m_0}, \psi) > 0$$

3. For each $\psi \in \Psi_m$, $y_1^n \in \mathbb{R}^n$ there exists $i = 1, \ldots I_{\varepsilon,m}$,

$$\left| \frac{\log p_{\psi_i}(Y_1^n|Y_0 = y_0) - \log p_\psi(Y_1^n|Y_0 = y_0)}{n} \right| \leq \varepsilon.$$

In the following theorem we prove that the estimator $\widehat{m}$ underestimates the number of states $m_0$.

**Theorem 4.2.** *Assume* (**S1**-**S6**) *and* (**H1**-**H3**). *Set* $\rho > 2$, *and for each* $n \geq 4$, $m \geq 1$

$$pen(n, m) = \sum_{l=1}^{m} \frac{l(l+1) + \rho}{2} \log n + \sum_{l=1}^{m} c_l(n) + \sum_{l=1}^{m} e_l(n) + m(m+1)\phi(n) \log n,$$

*where* $\phi(n) = o(n)$. *Then, for each* $m \leq m_0$ *it holds that* $\widehat{m}_m \leq m_0$ *a.s* $- \mathbb{P}_{\psi_0}$.

## 5. Proofs

**Proof of Theorem 3.1.**

We observe that

$$\int_\Psi p_\psi(Y_1^n|Y_0 = y_0)p(\psi)d\psi$$

$$= \sum_{x_1^n} \int_\Theta \int_\Sigma \int_\mathcal{P} p_{\theta,\sigma^2}(Y_1^n|Y_0 = y_0, x_1^n)p_A(x_1^n)p(A)p(\theta)p(\sigma^2)dAd\theta d\sigma^2$$

$$= \sum_{x_1^n} \int_\Sigma \int_\Theta p_\psi(Y_1^n|Y_0 = y_0, x_1^n)p(\theta)d\theta d\sigma^2 \int_\mathcal{P} p_A(x_1^n)p(A)dA$$

$$= \sum_{x_1^n} q_m(Y_1^n|Y_0 = y_0, x_1^n)q_m(x_1^n). \tag{5.1}$$

Hence, the Theorem can be proved by finding constants $C_1, C_2$ such that:

$$p_\theta(Y_1^n|Y_0 = y_0) \leq C_1 q_m(Y_1^n|Y_0 = y_0, x_1^n) \tag{5.2}$$

$$p_A(x_1^n) \leq C_2 q_m(x_1^n). \tag{5.3}$$

Thus, taking into account equations (5.1) and (3.1)

$$\begin{aligned} p_\psi(Y_1^n|Y_0 = y_0) &= \sum_{x_1^n} p_{\theta,\sigma^2}(Y_1^n|Y_0 = y_0, x_1^n)p_A(x_1^n) \\ &\leq C_1 C_2 \sum_{x_1^n} q_m(Y_1^n|x_1^n)q_m(x_1^n) \\ &= C_1 C_2 q_m(Y_1^n). \end{aligned}$$

Let us evaluate $q_m(x_1^n)$ following the proof given in the Appendix of Ref. [16]. Consider

$$q_m(x_1^n) = \prod_{i=1}^{m} \left[ \frac{\Gamma(m/2)}{\Gamma(n_i + 1/2)} \left( \prod_{i=1}^{m} \frac{\Gamma(n_{ij} + 1/2)}{\Gamma(1/2)} \right) \right]$$

and

$$\frac{p_A(x_1^n)}{q_m(x_1^n)} \leq \frac{\prod_{i=1}^{m} \prod_{j=1}^{m} \left(\frac{n_{ij}}{n_i}\right)^{n_{ij}}}{\prod_{i=1}^{m} \left[ \frac{\Gamma(m/2)}{\Gamma(n_i+1/2)} \left( \prod_{i=1}^{m} \frac{\Gamma(n_{ij}+1/2)}{\Gamma(1/2)} \right) \right]}. \tag{5.4}$$

The right-hand-side of equation (5.4) does not exceed

$$\left[ \frac{\Gamma(n + m/2)\Gamma(1/2)}{\Gamma(m/2)\Gamma(n + 1/2)} \right]^m.$$

Gassiat and Boucheron [10] showed that

$$m \log \left[ \frac{\Gamma(n+m/2)\Gamma(1/2)}{\Gamma(m/2)\Gamma(n+1/2)} \right] \leq \frac{m(m-1)}{2} \log n + c_m(n),$$

for $n \geq 4$, $c_m(n)$ one selects:

$$\log m - m \left( \log \frac{\Gamma(m/2)}{\Gamma(1/2)} - \frac{m(m-1)}{4n} + \frac{1}{12n} \right).$$

It follows:

$$\frac{p_A(x_1^n)}{q_m(x_1^n)} \leq n^{m(m-1)/2} e^{c_m(n)}. \tag{5.5}$$

What remains is to evaluate the quotient between $p_{\theta,\sigma^2}(Y_1^n | Y_0 = y_0, x_1^n, \theta, \sigma^2)$ and $q_m(Y_1^n | Y_0 = y_0, x_1^n)$. Let us start with the evaluation of $q_m$.

$$q_m(Y_1^n | Y_0 = y_0, x_1^n)$$
$$= \int \prod_{i=1}^{m} (2\pi\sigma_i^2)^{-n_i/2} e^{\left(-\frac{1}{2\sigma_i^2}(\mathbf{Y}_{I_i} - \mathbf{W}_i\theta_i)^t(\mathbf{Y}_{I_i} - \mathbf{W}_i\theta_i)\right)}$$
$$\times \frac{1}{2\pi\tau^2\sigma_i^2} e^{\left(-\frac{\theta_i^t\theta_i}{2\tau^2\sigma_i^2}\right)} \left(\frac{u_0}{2}\right)^{v_0/2} \frac{(\sigma_i^2)^{-(1+v_0/2)}}{\Gamma(v_0/2)} e^{\left(-\frac{u_0}{2\sigma_i^2}\right)} d\theta_i d\sigma_i^2.$$

As a result of the evaluation of the mixture, upon integration over the variables $\theta$ y $\sigma^2$, is:

$$q_m(Y_1^n | Y_0 = y_0, x_1^n)$$
$$= \prod_{i=1}^{m} \frac{\sqrt{\det(\mathbf{M}_i)}}{(2\pi)^{n_i}\tau^2} \left(\frac{u_0}{2}\right)^{v_0/2} \frac{2^{(v_0+n_i)/2}}{\Gamma(u_0/2)} (\mathbf{Y}_{I_i}^t P_i \mathbf{Y}_{I_i} + u_0)^{-(v_0+n_i)/2} \Gamma\left(\frac{v_0+n_i}{2}\right)$$

Now, setting $u_0 \to 0$ and $v_0 \to 0$ (which means that in the limit we consider *a priori* distributions which are not informative for $\sigma^2$ although they are improper)

$$q_m(Y_1^n | Y_0 = y_0, x_1^n) = \prod_{i=1}^{m} \frac{\sqrt{\det(\mathbf{M}_i)} 2^{n_i/2}}{(2\pi)^{n_i}\tau^2} \left(\mathbf{Y}_{I_i}^t P_i \mathbf{Y}_{I_i}\right)^{-n_i/2} \Gamma(n_i/2).$$

Again, introducing conditions with respect to $Y_1^n = y_1^n$ and $x_1^n$, as the model is both linear and Gaussian, the estimators ML for $1 \leq i \leq m$ are

$$\begin{aligned}\widehat{\theta}_i &= (\mathbf{W}_i^t \mathbf{W}_i)^{-1} \mathbf{W}_i^t \mathbf{Y}_{I_i} \\ \widehat{\sigma}_i^2 &= \frac{1}{n_i}(\mathbf{Y}_{I_i}^t \mathbf{Y}_{I_i} - \widehat{\theta}_i^t \mathbf{W}_i^t \mathbf{Y}_{I_i})\end{aligned}$$

Taking into account that $p_{\theta,\sigma^2}(Y_1^n|Y_0 = y_0, x_1^n) \leq p_{\widehat{\theta},\widehat{\sigma}^2}(Y_1^n|Y_0 = y_0, x_1^n)$ and that the right-hand-side of the inequality satisfies

$$\begin{aligned}p_{\widehat{\theta},\widehat{\sigma}^2}(Y_1^n|Y_0 = y_0, x_1^n) &= \prod_{i=1}^m (2\pi\widehat{\sigma}_i^2)^{-n_i/2} e^{-n_i/2} \\ &= \prod_{i=1}^m (2\pi)^{-n_i/2} e^{-n_i/2} n_i^{n_i/2} (\mathbf{Y}_{I_i}^t \mathbf{B}_i \mathbf{Y}_{I_i})^{-n_i/2}.\end{aligned}$$

We arrive at the following expression for the density quotient:

$$\frac{p_{\theta,\sigma^2}(Y_1^n|y_0, x_1^n)}{q_m(Y_1^n|Y_0 = y_0, x_1^n)} \leq \prod_{i=1}^m \frac{n_i^{n_i/2} \pi^{n_i/2}}{e^{n_i/2} \Gamma(n_i/2)} \left\{ \frac{\mathbf{Y}_{I_i}^T \mathbf{P}_i \mathbf{Y}_{I_i}}{\mathbf{Y}_{I_i}^t \mathbf{B}_i \mathbf{Y}_{I_i}} \right\}^{n_i/2} \tau^2 \sqrt{\det(\mathbf{M}_i^{-1})}.$$

Taking logarithms of both sides of the inequality, we have that:

$$\begin{aligned}\log \frac{p_{\theta,\sigma^2}(Y_1^n|Y_0 = y_0, x_1^n)}{q_m(Y_1^n|Y_0 = y_0, x_1^n)} &\leq \sum_{i=1}^m \log(d_i) + \sum_{i=1}^m \frac{n_i}{2} \log \frac{\mathbf{Y}_{I_i}^t \mathbf{P}_i \mathbf{Y}_{I_i}}{\mathbf{Y}_{I_i}^t \mathbf{B}_i \mathbf{Y}_{I_i}} \\ &\quad + \sum_{i=1}^m \log \tau^2 \sqrt{\det(\mathbf{M}_i^{-1})} \\ &= T_1 + T_2 + T_3 \quad \text{(say.)}\end{aligned}$$

Let us notice that the right-hand side of the former inequality satisfies the following bounds: For term $T_1$ we have

$$\sum_{i=1}^m \log\left(\frac{n_i^{n_i/2} \pi^{n_i/2}}{e^{n_i/2} \Gamma(n_i/2)}\right) \leq \frac{n}{2} + \frac{1}{2}\log\left(\frac{n}{2}\right) - \frac{m \log(2\pi)}{2}.$$

For term $T_2$

$$\sum_{i=1}^m \frac{n_i}{2} \log \frac{\mathbf{Y}_{I_i}^t \mathbf{P}_i \mathbf{Y}_{I_i}}{\mathbf{Y}_{I_i}^t \mathbf{B}_i \mathbf{Y}_{I_i}} \leq \frac{nm}{2} \log \frac{\mathbf{Y}_{I_K}^t \mathbf{P}_k \mathbf{Y}_{I_k}}{\mathbf{Y}_{I_k}^t \mathbf{B}_k \mathbf{Y}_{I_k}},$$

and for term $T_3$

$$\tau^4 \det(\mathbf{M}_i^{-1}) = 1 + \tau^4 n_i \sum_{k \in I_i} Y_{k-1}^2 - \tau^4 \left(\sum_{k \in I_i} Y_{k-1}\right)^2 + \tau^2 + \tau^2 \sum_{k \in I_i} Y_{k-1},$$

we write the first term of the inequality

$$\sum_{i=1}^{m} \log \tau^2 \sqrt{\det(\mathbf{M}_i^{-1})}$$

$$= \sum_{i=1}^{m} \log \sqrt{1 + \tau^4 n_i \sum_{k \in I_i} Y_{k-1}^2 - \tau^4 \left(\sum_{k \in I_i} Y_{k-1}\right)^2 + \tau^2 + \tau^2 \sum_{k \in I_i} Y_{k-1}}$$

$$= \sum_{i=1}^{m} \log \sqrt{1 + V_i} \quad \text{(say.)}$$

Making use both of convexity and the Ergodic Theorem we see that the following relation is satisfied:

$$\sum_{i=1}^{m} \log \sqrt{1 + V_i} \leq \log \left(1 + \frac{1}{m} \sum_{i=1}^{m} V_i\right)^{m/2} = \log \left(1 + \frac{\tau^4 n^2}{m} \sum_{i=1}^{m} (\lambda_i \sigma_i)^2\right)^{m/2} \quad a.s.$$

Substituting the calculated bounds

$$\log \frac{p_{\theta,\sigma^2}(Y_1^n | Y_0 = y_0, x_1^n)}{q_m(Y_1^n | Y_0 = y_0, x_1^n)}$$

$$\leq -\frac{m}{2}(\log(2) + \log(2\pi) + \log(n)) + \frac{\log 2\pi}{2} n + \frac{n}{2} \log \frac{\mathbf{Y}_{I_k}^t \mathbf{P}_k \mathbf{Y}_{I_k}}{\mathbf{Y}_{I_k}^t \mathbf{B}_k \mathbf{Y}_{I_k}}$$

$$+ \log \left(1 + \frac{\tau^4 n^2}{m} \sum_{i=1}^{m} (\lambda_i \sigma_i)^2\right)^{m/2}$$

□

**Proof of Lemma 4.1.**

We work directly with the extended Markov chain $\{(Y_n, X_n)\}$. Let $h(\psi) = \frac{1}{n} \log p_\psi(Y_0^n, x_1^n)$ and let $\psi, \psi' \in \Psi$. We prove that for each $\varepsilon > 0$ there exists a $\delta(\varepsilon) > 0$ such that:

$$\forall n \quad |h_n(\psi) - h_n(\psi')| \leq \varepsilon \quad \text{si} \quad \|\psi - \psi'\| < \delta(\varepsilon).$$

Complete likelihood is written as

$$p_\psi(Y_0^n, x_1^n)$$
$$= \prod_{k=1}^{n} \prod_{i,j=1}^{m} a_{ij}^{\mathbb{1}_{i,j}(x_k, x_{k+1})} \prod_{i=1}^{m} \frac{1}{(2\pi\sigma_i^2)^{n_i/2}} e^{\left(-\frac{1}{2\sigma_i^2}(\mathbf{Y}_{I_i} - \mathbf{W}_i \theta_i)^t (\mathbf{Y}_{I_i} - \mathbf{W}_i \theta_i)\right)}$$

from where it ensures that

$$\begin{aligned}
|h_n(\psi) - h_n(\widetilde{\psi})| \\
\leq\ & \frac{1}{n}\sum_{i,j=1}^{m} n_{ij}|\log a_{ij} - \log \widetilde{a}_{ij}| + \frac{1}{2n}\sum_{i=1}^{m} n_i \left|\log \sigma_i^2 - \log \widetilde{\sigma_i^2}\right| \\
+\ & \frac{1}{n}\left|\sum_{i=1}^{m}\left(\frac{1}{2\sigma_i^2} - \frac{1}{2\widetilde{\sigma_i^2}}\right)\mathbf{Y}_{I_i}^t \mathbf{Y}_{I_i}\right| + \frac{1}{n}\left|\sum_{i=1}^{m} \mathbf{Y}_{I_i}^t \mathbf{W}_i\left(\frac{\theta_i}{\sigma_i^2} - \frac{\widetilde{\theta_i}}{\widetilde{\sigma_i^2}}\right)\right| \\
+\ & \frac{1}{n}\left|\sum_{i=1}^{m}\left(\frac{\theta_i}{\sigma_i^2} - \frac{\widetilde{\theta_i}}{\widetilde{\sigma_i^2}}\right)^t \mathbf{W}_i^t \mathbf{W}_i \left(\frac{\theta_i}{\sigma_i^2} - \frac{\widetilde{\theta_i}}{\widetilde{\sigma_i^2}}\right)\right| \\
=\ & T_1 + T_2 + T_3 + T_4 + T_5 \qquad \text{(say.)}
\end{aligned} \tag{5.6}$$

The right-hand-side of the inequality (5.6) can be bounded in the following way

- Since $n_{ij}/n \leq 1$, $n_i/n \leq 1$ and the parameters $a_{ij}, \widetilde{a}_{ij}, \sigma_i^2, \widetilde{\sigma_i^2}$ are lower bounded (for **S6**), there exist a constant $C_1$ such that the terms $T_1$ and $T_2$ of (5.6) are upper bounded by $C_1\|\psi - \widetilde{\psi}\|$.
- Due to compactness of the parameters space (**S6**), there exist a constant $C_2$ such that term $T_3$ of (5.6) are upper bounded by $C_2\|\sigma^2 - \widetilde{\sigma^2}\|\frac{1}{n}\sum_{k=1}^{n} Y_k$. The stability condition (**S4**), and the existence of the moments of $e_1$ (**S3**) implies (see Yao and Atalli [22]), by the Ergodic Theorem, that the terms of the $1/n \sum_{k=1}^{n} g(Y_k)$ are controlled. Hence

$$C_2\|\sigma^2 - \widetilde{\sigma^2}\|\frac{1}{n}\sum_{k=1}^{n} Y_k \leq C_3\|\psi - \widetilde{\psi}\| \quad \text{a.s.}$$

- By the same argument of compactness (**S5** and **S6**) we have

$$T_4 \leq C_4\|\psi - \widetilde{\psi}\|\left|\frac{1}{n}\sum_{k=1}^{n} Y_k + \frac{1}{n}\sum_{i=1}^{m}\sum_{k\in I_i}^{n} Y_k Y_{k-1}\right|$$

and again, following the Ergodic Theorem, the right side of the above inequality is upper bounded by $C_4\|\psi - \widetilde{\psi}\|$ a.s.
- By the Cauchy-Schwarz-Bunyakowski inequality

$$\begin{aligned}
T_5 &\leq \frac{1}{n}\sum_{i=1}^{m}\left\|\left(\frac{\theta_i}{\sigma_i^2} - \frac{\widetilde{\theta_i}}{\widetilde{\sigma_i^2}}\right)\right\| \|\mathbf{W}_i^t \mathbf{W}_i\| \\
&\leq C_5\|\psi - \widetilde{\psi}\|\frac{1}{n}\sum_{i=1}^{m}\|\mathbf{W}_i^t \mathbf{W}_i\|.
\end{aligned}$$

Now, the norm of the symmetric matrix $\mathbf{W}_i^t \mathbf{W}_i$ is given by the absolute value of of the largest real eigenvalue, which in the present case is

$$\frac{\mathrm{tr}(\mathbf{W}_i^t \mathbf{W}_i) + \sqrt{\mathrm{tr}(\mathbf{W}_i^t \mathbf{W}_i)^2 - 4\det \mathbf{W}_i^t \mathbf{W}_i}}{2}.$$

Since $\det \mathbf{W}_i^t \mathbf{W}_i$ is positive,

$$\frac{\operatorname{tr}(\mathbf{W}_i^t \mathbf{W}_i) + \sqrt{\operatorname{tr}(\mathbf{W}_i^t \mathbf{W}_i)^2 - 4 \det \mathbf{W}_i^t \mathbf{W}_i}}{2} \leq \operatorname{tr}(\mathbf{W}_i^t \mathbf{W}_i).$$

Since $\operatorname{tr}(\mathbf{W}_i^t \mathbf{W}_i) = n_i + \sum_{k \in I_i} Y_k^2$, then

$$\frac{1}{n} \sum_{i=1}^{m} \|\mathbf{W}_i^t \mathbf{W}_i\| \leq 1 + \frac{1}{n} \sum_{k=1}^{n} Y_k^2.$$

Thus, the last term of (5.6) is smaller than $C_5 \|\psi - \widetilde{\psi}\|$.

We thus reach the conclusions that there exists a constant $C$ such that

$$|h_n(\psi) - h_n(\psi')| \leq C \|\psi - \psi'\|, \quad \text{a.s.}$$

this implies that $h_n$ is an equicontinuous series. In order to return to $\{Y_n\}$ we note that

$$\left| \frac{1}{n} \log \frac{p_\psi(Y_0^n, x_1^n)}{p_{\psi'}(Y_0^n, x_1^n)} \right| \leq \varepsilon.$$

from where we have

$$p_{\psi'}(Y_0^n, x_1^n) \leq e^{(\varepsilon n)} p_\psi(Y_0^n, x_1^n).$$

and then, adding over $x_1^n$

$$p_{\psi'}(Y_1^n | Y_0 = y_0) = \sum_{x_1^n} p_{\psi'}(Y_0^n, x_1^n) \leq e^{(\varepsilon n)} \sum_{x_1^n} p_\psi(Y_0^n, x_1^n) = p_\psi(Y_1^n | Y_0 = y_0)$$

from where it follows

$$\left| \frac{1}{n} \log \frac{p_{\psi'}(Y_1^n | Y_0 = y_0)}{p_\psi(Y_1^n | Y_0 = y_0)} \right| \leq \varepsilon.$$

□

**Proof of Lemma 4.2.**

The first part follows from proposition 2.9 of [15].

To prove the second part, we follow Leroux lemma (see [3], Lemma 8, p. 21), for every $\psi \in \Psi_{m_0}$ such that $p_\psi \neq p_{\psi_{m_0}}$, there exists a neighborhood $O_\psi$ and $\varepsilon > 0$ tal que $\inf_{\psi \in O_\psi} K(\psi_{m_o}, \psi) > \varepsilon$. Since, however, $\Psi_{m_0-1}$ is compact, it is subcovering by a finite union $O_{\psi_1}, \ldots, O_{\psi_I}$ (each one of them is associated to a $\varepsilon_i > 0$); hence,

$$\inf_{\psi \in \Psi_{m_0-1}} K(\psi, \psi_0) \geq \min_{i \leq I} \inf_{\psi \in O_{\psi_i}} K(\psi, \psi_0) \geq \min_{i \leq I} \varepsilon_i > 0.$$

In order to carry our the third part of this proof, let $\{B_\varepsilon(\psi) : \psi \in \Psi_m\}$ be a covering of $\Psi_m$ by open balls. Due to the compactness of $\Psi_m$ there exists a finite subcovering $B_\varepsilon(\psi_1), \ldots, B_\varepsilon(\psi_I)$. Thus, for every $\psi \in \Psi_m$ there exists $i \in \{1, \ldots, I\}$ such that from Lemma 1.1,

$$\left| \frac{\log p_{\psi_i}(Y_1^n | Y_0 = y_0) - \log p_\psi(Y_1^n | Y_0 = y_0)}{n} \right| \leq \varepsilon.$$

□

**Proof of Theorem 4.1.**

We use the fact that $\mathbb{P}(\widehat{m}(n) < m_0 \text{ i.o}) \leq \sum_{m=1}^{m_0-1} \mathbb{P}(\widehat{m}(n) = m)$. We prove that $\mathbb{P}(\widehat{m}(n) = m) = 0$. Indeed,

$$\begin{aligned}
\mathbb{P}(\widehat{m}(n) = m) \\
\leq \ & \mathbb{P}\Big( \sup_{\psi \in \Psi_m} \log p_\psi - pen(n,m) \geq \log p_{\psi_{m_0}} - pen(n, m_0) \Big) \\
\leq \ & \mathbb{P}\Big( \sup_{\psi \in \Psi_m} \log p_\psi \geq \log p_{\psi_{m_0}} - pen(n, m_0) + pen(n, m) \Big), \quad (5.7)
\end{aligned}$$

since $\psi \in \Psi_m$ according to Lemma 4.2 there exists $1 \leq i \leq I$ such that $\log p_{\psi_m} < n\varepsilon + \log p_{\psi_i}$; hence it follows from the (5.7) that

$$\begin{aligned}
\mathbb{P}(\widehat{m}(n) = m) &\leq \mathbb{P}\Big( \max_{i \leq I} \log p_{\psi_i} \geq \log p_{\psi_{m_0}} - pen(n, m_0) - n\varepsilon \Big) \\
&\leq \sum_{i=1}^{I} \mathbb{P}\left( \frac{\log p_{\psi_i} - \log p_{\psi_{m_0}}}{n} \geq -\frac{pen(n, m_0)}{n} - \varepsilon \right)
\end{aligned}$$

and again, according to Lemma 4.2,

$$\lim_{n \to \infty} \frac{\log p_{\psi_i} - \log p_{\psi_{m_0}}}{n} = -K(\psi_i, \psi_0)$$

and by hypothesis $\lim_{n \to \infty} \frac{pen(n,m_0)}{n} = 0$ from where it follows that,

$$\mathbb{P}(\widehat{m}(n) = m) \leq \sum_{i=1}^{I} \mathbb{P}\left( \varepsilon < K(\psi_i, \psi_0) \leq \varepsilon \right) = 0.$$

□

**Proof of Theorem 4.2.**

Let us define the set

$$A_n = \left\{ \frac{\mathbf{Y}_{I_k}^t \mathbf{P}_k \mathbf{Y}_{I_k}}{\mathbf{Y}_{I_k}^t \mathbf{B}_k \mathbf{Y}_{I_k}} \leq t_n \right\}$$

and

$$\Delta_{n,m} = c_m(n) + d_m(n) + e_m(n) + \frac{nm}{2} \log \frac{\mathbf{Y}_k^t \mathbf{P}_k \mathbf{Y}_k}{\mathbf{Y}_k^t \mathbf{B}_k \mathbf{Y}_k} + pen(n, m_0) - pen(n, m).$$

We note that

$$\mathbb{P}_{\psi_0}(\widehat{m}(n) > m_0) \leq \sum_{m > m_0} \mathbb{P}_{\psi_0}(\widehat{m}(n) > m_0, A_n) + \mathbb{P}_{\psi_0}(A_n^c)$$

and

$$\mathbb{P}_{\psi_{m_0}}(\widehat{m} = m, A_n)$$

$$\stackrel{(a)}{\leq} \mathbb{P}_{\psi_{m_0}}(\log p_{\psi_{m_0}}(Y_1^n | Y_0 = y_0) \leq \sup_{\psi \in \Psi_m} \log p_{\psi_m}(Y_1^n | Y_0 = y_0)$$
$$+ pen(n, m_0) - pen(n, m), A_n)$$

$$\stackrel{(b)}{\leq} \mathbb{P}_{\psi_{m_0}}\left(\log \frac{p_{\psi_{m_0}}(Y_1^n | Y_0 = y_0)}{q_m(y_1^n)} \leq \Delta_{m,n}, A_n\right)$$

$$= \int_{Y_1^n} \mathbb{1}\left(\log \frac{p_{\psi_{m_0}}(Y_1^n = y_1^n | y_0)}{q_m(Y_1^n = y_1^n)} \leq \Delta_{m,n}, A_n\right)$$

$$\times \frac{p_{\psi_{m_0}}(Y_1^n = y_1^n | Y_0 = y_0)}{q_m(Y_1^n = y_1^n)} q_m(Y_1^n = y_1^n) dy_1^n$$

$$\leq \exp\left(\frac{m(m+1)}{2}\log(n) + c_m(n) + d(n) + e_m(n)\right.$$
$$\left. + \frac{n\log(t_n)}{2} + pen(n, m_0) - pen(n, m)\right)$$

where (a) is a consequence of the way the PML estimator is defined (4.2) and (b), of Theorem 3.1.

In what follows we consider the coefficient $\mathbf{Y}_{I_k}^t \mathbf{P}_k \mathbf{Y}_{I_k} / \mathbf{Y}_{I_k}^t \mathbf{B}_k \mathbf{Y}_{I_k}$. Conditions with respect to $Y_1^n = y_1^n$ and $x_1^n$, as the model is both linear and Gaussian (3.2) then $\mathbf{Y}_{I_k}^t \mathbf{P}_k \mathbf{Y}_{I_k}$ has a $\chi^2(n_k, \gamma)$ distribution, where $\gamma = (1/2)\theta_k^t \mathbf{W}^t \mathbf{P}_k \mathbf{W} \theta_k$ is the non-centrality parameter; further, we assume that $\mathbf{P}_k$ has a maximum rank. Moreover, $\chi^2(n_k, 1/2\theta_k^t \mathbf{W}^t \mathbf{P}_k \mathbf{W} \theta_k)$ can be approximated by a $\chi_r^2$ having the same mean and the same variance with $r = (n_k + 2\gamma)^2/(n_k + 4\gamma)$. For the denominator, if assume $\mathbf{B}_k$ to have full range, then $\mathbf{Y}_{I_k}^t \mathbf{B}_k \mathbf{Y}_{I_k}$ distributes $\chi_{n_k}^2$ (see Searle [20],§2, págs. 49–53).

On the other hand,

$$\gamma = \frac{1}{2}\tau^{-2}\mathbf{Y}_{I_k}^t \mathbf{W}_k(\mathbf{W}_k^t \mathbf{W}_k)^{-1}\mathbf{M}_k \mathbf{W}_k^t \mathbf{Y}_{I_k} \approx \frac{1}{2}\tau^{-2}\|\mathbf{Y}_{I_k}^t \mathbf{W}_k\|^2,$$

substituting in $r$, we have:

$$r = \frac{(n_k + 2\gamma)^2}{(n_k + 4\gamma)} = \frac{(n_k + \tau^{-2}\|\mathbf{Y}_{I_k}^t \mathbf{W}_k\|^2)^2}{(n_k + 2\tau^{-2}\|\mathbf{Y}_{I_k}^t \mathbf{W}_k\|^2)} = o\left(\frac{n^2}{2\tau^2}\right) \text{ a.s.}$$

Since we notice that $n_k/n \to \lambda_k$ a.s (**S1**),

$$\mathbb{P}\left(\frac{\mathbf{Y}_{I_k}^T t\mathbf{P}_k \mathbf{Y}_{I_k}}{\mathbf{Y}_{I_k}^t \mathbf{B}_k \mathbf{Y}_{I_k}} \geq t_n\right) \approx \mathbb{P}\left(\mathcal{F}_{n_k,r} \leq \frac{1}{t_n}\right)$$

$$= \frac{\Gamma\left(\frac{n_k+r}{2}\right)}{\Gamma\left(\frac{n_k}{2}\right)\Gamma\left(\frac{r}{2}\right)} \int_0^{1/t_n} u^{n_k/2-1}(1+u)^{-(n_k+r)/2} du$$

$$= \frac{\Gamma\left(\frac{n_k+r}{2}\right)}{\Gamma\left(\frac{n_k}{2}\right)\Gamma\left(\frac{r}{2}\right)} \frac{2}{n_k t_n^{\frac{n_k}{2}}}$$

$$\leq \frac{1}{\sqrt{\pi} n^{\frac{3n^2}{4\tau^2}} \lambda_k t_n^{\frac{n\lambda_k}{2}}}$$

and choosing $t_n = n^{\frac{3}{2\lambda_k \tau^2}}$ we have that

$$\mathbb{P}\left(\frac{\mathbf{Y}_{I_k}^T t\mathbf{P}_k \mathbf{Y}_{I_k}}{\mathbf{Y}_{I_k}^t \mathbf{B}_k \mathbf{Y}_{I_k}} \geq t_n\right) \leq \frac{1}{\sqrt{\pi}\lambda_k n^{\frac{3(n^2+n)}{4\tau^2}}}$$

We have proven that $\mathbf{Y}_{I_k}^t \mathbf{P}_k \mathbf{Y}_{I_k}/\mathbf{Y}_{I_k}^t \mathbf{B}_k \mathbf{Y}_{I_k}$ is bounded in probability with a rate $n^{\frac{3}{2\lambda_k \tau^2}}$. We still have to determine the bounds for $\Delta_{nm}$.

Using the definition of function $pen(n,m)$ we have

$$\Delta_{nm} \leq -\frac{\rho}{2}(m-m_0)\log(n) - \sum_{l=m_0+1}^{m-1} \frac{l(l+1)}{2}\log n - \sum_{l=m_0+1}^{m-1} c_l(n)$$

$$- \sum_{l=m_0+1}^{m-1} e_l(n) + mn\log n + \frac{3m\phi(n)\log(n)}{4\lambda_k \tau^2}$$

$$+ \frac{m_0(m_0+1)}{2}\phi(n)\log n - \frac{m(m+1)}{2}\phi(n)\log n$$

for $m = m_0 + 1$ we have that

$$-\sum_{l=m_0+1}^{m-1} \frac{l(l+1)}{2}\log n - \sum_{l=m_0+1}^{m-1} c_l(n) - \sum_{l=m_0 1}^{m-1} e_l(n) = 0.$$

We select $\tau^2 = \frac{3}{4\lambda_k}$ as:

$$\Delta_{nm} \leq -\frac{\rho}{2}(m-m_0)\log(n) + \left[(m_0-m)(m_0+m) + m_0 + \frac{3m}{4\lambda_k \tau^2}\right]\phi(n)\log n$$

$$\leq -\frac{\rho}{2}(m-m_0)\log(n)$$

Therefore:

$$\mathbb{P}_{\psi_{m_0}}\left(\widehat{m} = m,\ A_n\right) \leq e^{\left(-\frac{\rho}{2}(m-m_0)\log(n)\right)} = O(n^{-\rho/2}),$$

and $\mathbb{P}_{\psi_0}(A_n^c) = O(e^{-n\log n})$ hence $\mathbb{P}_{\psi_0}(\widehat{m}(n) > m_0) = O(n^{-\rho/2} + e^{-n\log n})$ thus, in view of Borel-Cantelli's Lemma $\widehat{m}(n) \leq m_0$ a.s. □